\documentclass{amsart}
\usepackage{latexsym,amsmath,amssymb}

\theoremstyle{definition}
\theoremstyle{remark}

\numberwithin{equation}{section}

\begin{document}

\title[Weighted conditional expectation type operators]
{Decompositions of weighted conditional expectation type
operators}

\author{\sc\bf Y. Estaremi }
\address{\sc Y. Estaremi }
\email{estaremi@gmail.com - yestaremi@pnu.ac.ir}

\address{Department of Mathematics, University of Payame noor, p. o. box: 19395-3697, Tehran, Iran.}

\thanks{}

\thanks{}



\date{}

\dedicatory{}

\commby{}

\begin{abstract}
In this paper we investigate boundedness, polar decomposition and
spectral decomposition of weighted conditional expectation type
operators on $L^2(\Sigma)$.

\noindent {}
\end{abstract}

\maketitle

\section{ \sc\bf Introduction and Preliminaries}

Let $(X,\Sigma,\mu)$ be a complete $\sigma$-finite measure space.
For any sub-$\sigma$-finite algebra $\mathcal{A}\subseteq
 \Sigma$ with $1\leq p\leq \infty$, the $L^p$-space
$L^p(X,\mathcal{A},\mu_{\mid_{\mathcal{A}}})$ is abbreviated  by
$L^p(\mathcal{A})$, and its norm is denoted by $\|.\|_p$. All
comparisons between two functions or two sets are to be
interpreted as holding up to a $\mu$-null set. The support of a
measurable function $f$ is defined as $S(f)=\{x\in X; f(x)\neq
0\}$. We denote the vector space of all equivalence classes of
almost everywhere finite valued measurable functions on $X$ by
$L^0(\Sigma)$.

\vspace*{0.3cm} For a sub-$\sigma$-finite algebra
$\mathcal{A}\subseteq\Sigma$, the conditional expectation operator
associated with $\mathcal{A}$ is the mapping $f\rightarrow
E^{\mathcal{A}}f$, defined for all non-negative, measurable
function $f$ as well as for all $f\in L^p(\Sigma)$, $1\leq p\leq
\infty$, where $E^{\mathcal{A}}f$, by the Radon-Nikodym theorem,
is the unique $\mathcal{A}$-measurable function satisfying
$$\int_{A}fd\mu=\int_{A}E^{\mathcal{A}}fd\mu, \ \ \ \forall A\in \mathcal{A} .$$
As an operator on $L^{p}({\Sigma})$, $E^{\mathcal{A}}$ is
idempotent and $E^{\mathcal{A}}(L^p(\Sigma))=L^p(\mathcal{A})$. If
there is no possibility of confusion, we write $E(f)$ in place of
$E^{\mathcal{A}}(f)$. This operator will play a major role in our
work and we list here some of its useful properties:

\vspace*{0.2cm} \noindent $\bullet$ \  If $g$ is
$\mathcal{A}$-measurable, then $E(fg)=E(f)g$.

\noindent $\bullet$ \ $|E(f)|^p\leq E(|f|^p)$.

\noindent $\bullet$ \ If $f\geq 0$, then $E(f)\geq 0$; if $f>0$,
then $E(f)>0$.

\noindent $\bullet$ \ $|E(fg)|\leq
E(|f|^p)|^{\frac{1}{p}}E(|g|^{q})|^{\frac{1}{q}}$, where
$\frac{1}{p}+\frac{1}{q}=1$ (H\"{o}lder inequality).

\noindent $\bullet$ \ For each $f\geq 0$, $S(f)\subseteq S(E(f))$.

\vspace*{0.2cm}\noindent A detailed discussion and verification of
most of these properties may be found in \cite{rao}. We recall
that an $\mathcal{A}$-atom of the measure $\mu$ is an element
$A\in\mathcal{A}$ with $\mu(A)>0$ such that for each
$F\in\mathcal{A}$, if $F\subseteq A$, then either $\mu(F)=0$ or
$\mu(F)=\mu(A)$. A measure space $(X,\Sigma,\mu)$ with no atoms is
called a non-atomic measure space. It is well-known fact that
every $\sigma$-finite measure space $(X,
\mathcal{A},\mu_{\mid_{\mathcal{A}}})$ can be partitioned uniquely
as $X=\left (\bigcup_{n\in\mathbb{N}}A_n\right )\cup B$, where
$\{A_n\}_{n\in\mathbb{N}}$ is a countable collection of pairwise
disjoint $\mathcal{A}$-atoms and $B$, being disjoint from each
$A_n$, is non-atomic (see \cite{z}).

\vspace*{0.2cm} Combinations of conditional expectation operators
and multiplication operators  appear often in the study of other
operators such as multiplication operators, weighted composition
operators and integral operators. Specifically, in \cite{mo},
S.-T. C. Moy characterized all operators on $L^p$ of the form
$f\rightarrow E(fg)$ for $g$ in $L^q$ with $E(|g|)$ bounded.
Eleven years later, R. G. Douglas, \cite{dou}, analyzed positive
projections on $L^{1}$ and many of his characterizations are in
terms of combinations of multiplications and conditional
expectations. More recently, P.G. Dodds, C.B. Huijsmans and B. De
Pagter, \cite{dhd}, extended these characterizations to the
setting of function ideals and vector lattices. J. Herron
presented some assertions about the operator $EM_u$ on $L^p$
spaces in \cite{he, her}. \\
In \cite{e, ej} we investigated some classic properties of
multiplication conditional expectation operators $M_wEM_u$ on
$L^p$ spaces. Let $f\in L^0(\Sigma)$, then $f$ is said to be
conditionable with respect to $E$ if $f\in\mathcal{D}(E):=\{g\in
L^0(\Sigma): E(|g|)\in L^0(\mathcal{A})\}$. Throughout this paper
we take $u$ and $w$ in $\mathcal{D}(E)$.  In this paper we present
some results on the boundedness, polar decomposition and spectral
decomposition of this operators in $L^2(\Sigma)$, using different
methods than those employed in \cite{ej}.

\section{ \sc\bf Polar decomposition }

\vspace*{0.3cm} {\bf Theorem 2.1.} The operator $T=M_{w}EM_{u}:
L^2(\Sigma)\rightarrow L^2(\Sigma)$ is bounded if and only if
$(E(|w|^2)^{\frac{1}{2}}) (E(|u|^2)^{\frac{1}{2}})\in
L^{\infty}(\mathcal{A})$ and in this case
$\|T\|=\|(E(|w|^2)^{\frac{1}{2}})
(E(|u|^2)^{\frac{1}{2}})\|_{\infty}$.

\vspace*{0.3cm} {\bf Proof}
 Suppose that $(E(|w|^2)^{\frac{1}{2}})
(E(|u|^2)^{\frac{1}{2}})\in L^{\infty}(A)$. Let $f\in
L^2(\Sigma)$. Then

$$\|T(f)\|^{2}_{2}=\int_{X}|wE(uf)|^2d\mu=\int_{X}E(|w|^2)|E(uf)|^2d\mu\leq\int_{X}E(|w|^2)E(|u|^2)E(|f|^2)|d\mu.$$

Since $|E(uf)|\leq(E(|u|^2)^{\frac{1}{2}}(E(|f|^2)^{\frac{1}{2}}$.
Thus

$$\|T\|\leq
\|E(|w|^2)^{\frac{1}{2}}(E(|u|^2))^{\frac{1}{2}}\|_{\infty}.$$\\

To prove the converse, let $T$ be bounded on $L^2(\Sigma)$ and
consider the case that $\mu(X)<\infty$.
 Then for
all $f\in L^2(\Sigma)$ we have

$$\|T(f)\|^{2}_{2}=\int_{X}|wE(uf)|^2d\mu=\int_{X}E(|w|^2)|E(uf)|^2d\mu$$$$\leq\|T\|^2\int_{X}|f|^2d\mu.$$

For each $n\in \mathbb{N}$, define $$E_{n}=\{x\in
X:|u(x)|(E(|w|^2))^{\frac{1}{2}}(x)\leq n\}.$$ Each $E_n$ is
$\Sigma$-measurable and $E_n\uparrow X$. Define $G_n=E_n\cap S$
for each $n\in \mathbb{N}$, where
$S=S(|u|(E(|w|^2))^{\frac{1}{2}})$. Let $A\in \mathcal{A}$ and
define $$f_n=\bar{u}(E(|w|^2))^{\frac{1}{2}} \chi_{G_n\cap A}$$
for each $n\in \mathbb{N}$. It is clear that $f_n\in
L^{\infty}(\Sigma)$ for all $n$ (which in our case implies $f_n\in
L^2(\Sigma)$). For each $n$,
\begin{align*}
\|T(f_n)\|^2_2&=\int_{X}E(|w|^2)|E(uf_n)|^2d\mu\\
&=\int_{X}(E(|w|^2))^2(E(|u|^2\chi_{G_n}.\chi_{A}))^2d\mu\\
&=\int_{A}[E(|w|^2)E(|u|^2\chi_{G_n})]^2d\mu\\
&\leq \|T\|^2\int_{X}|f_n|^2d\mu\\
&=\|T\|^2\int_{A}E(|u|^2\chi_{G_n})E(|w|^2)d\mu.
\end{align*}
Since $A$ is an arbitrary $\mathcal{A}$-measurable set and the
integrands are both $\mathcal{A}$-measurable functions, we have
$$[E(|w|^2)E(|u|^2\chi_{G_n})]^2\leq \|T\|^2E(|u|^2\chi_{G_n})E(|w|^2)$$ almost
everywhere. That is
$$[E(((E(|w|^2))^{\frac{1}{2}}|u|\chi_{E_n})^2\chi_{S})]^2\leq \|T\|^2E((|u|\chi_{E_n}(E(|w|^2))^{\frac{1}{2}})^q\chi_{S}).$$

Since
$$S=\sigma(|u|(E(|w|^2))^{\frac{1}{2}})=\sigma|u|^2E(|w|^2)$$
and
$$|u|^2E(|w|^2)\chi_{S}=|u|^2E(|w|^2),$$
we have
$$E(((E(|w|^2))^{\frac{1}{2}}
|u|\chi_{E_n})^2\chi_{S})\leq \|T\|^2.$$

Thus
$$(E|w|^2)^{\frac{1}{2}}(E|u|^2\chi_{E_n})^{\frac{1}{2}}\leq \|T\|.$$
This implies that
$(E|w|^2)^{\frac{1}{2}}(E|u|^2\chi_{E_n})^{\frac{1}{2}}\in
L^{\infty}(\mathcal{A})$ and
$$\|(E|w|^2)^{\frac{1}{2}}(E|u|^2)^{\frac{1}{2}}\|_{\infty}\leq
\|T\|.$$ Moreover, since  $E_m\uparrow X$, the conditional
expectation version of the monotone convergence theorem implies
$\|(E(|w|^2))^{\frac{1}{2}}
(E(|u|^2))^{\frac{1}{2}}\|_{\infty}\leq \|T\|$. $\hfill\Box$

\vspace*{0.3cm} {\bf proposition 2.2.} Let $g\in L^{\infty}(
\mathcal{A})$ and let $T=M_wEM_u:L^{2}(\Sigma)\rightarrow
L^{2}(\Sigma)$ be bounded. If $M_gT=0$, then $g=0$ on
$\sigma(E(|w|^{2})E(|u|^{2}))$.

\vspace*{0.3cm} {\bf Proof.} Let $f\in L^{2}(\Sigma)$. Then
$gwE(uf)=M_gT(f)=0$. Now, by Theorem 2.1
$$0=\|M_gT\|^{2}=\||g|^{2}E(|w|^{2})E(|u|^{2})\|{\infty},$$ which
implies that $|g|^{2}E(|w|^{2})E(|u|^{2})=0$, and so $g=0$ on
$\sigma(E(|w|^{2})E(|u|^{2}))$. $\hfill\Box$

\vspace*{0.3cm} {\bf Theorem 2.3.} The bounded operator
$T=M_wEM_u$ is a partial isometry if and only if
$E(|w|^{2})E(|u|^{2})=\chi_{A}$ for some $A\in\mathcal{A}$.

\vspace*{0.3cm} {\bf Proof.} Suppose $T$ is partial isometry. Then
$TT^{\ast}T=T$, that is $Tf=E(|w|^{2})E(|u|^{2})Tf$, and hence
$(E(|w|^{2})E(|u|^{2})-1)Tf=0$ for all $f\in L^{2}(\Sigma)$. Put
$S=S(E(|u|^2))$ and $G=S(E(|w|^2))$. By Proposition 2.2. we get
that $E(|w|^{2})E(|u|^{2})=1$ on $S\cap G$, which implies that
$E(|w|^{2})E(|u|^{2})=\chi_{A}$, where $A=S\cap G$.

\vspace*{0.3cm} Conversely, suppose that
$E(|w|^{2})E(|u|^{2})=\chi_{A}$ for some $A\in\mathcal{A}$. It
follows that $A=S\cap G$, and we have
$$TT^{\ast}T(f)=E(|w|^{2})E(|u|^{2})Tf=\chi_{S\cap G}wE(uf)=wE(uf),$$
where we have used the fact that $S(Tf)=S(|Tf|^{2})\subseteq S\cap
G$, which this is a consequence of H\"{o}lder's inequality for
conditional expectation $E$.$\hfill\Box$

 The spectrum of an
operator $A$ is the set
$$\sigma(A)=\{\lambda\in \mathbb{C}:A-\lambda I \ \ \  is \  not \
invertible\}.$$ It is well known that any bounded operator $A$ on
a Hilbert space $\mathcal{H}$ can be expressed in terms of its
polar decomposition: $A = UP$, where $U$ is a partial isometry and
$ P$ is a positive operator. (An operator is positive if $\langle
Pf,f\rangle \geq0$, for all $f\in \mathcal{H}$.) This
representation is unique under the condition that $ker U=ker P=ker
A$. Moreover, $P=|A|=(A^*A)^{\frac{1}{2}}$.

Let $q(z)$ be a polynomial with complex coefficients:
$q(z)=\sum^{N}_{n=0}\alpha_nz^n$. If $T$ is a bounded operator on
$L^2(\Sigma)$, then the operator $q(T)$ is defined by
$q(z)=\alpha_0I+\sum^{N}_{n=1}\alpha_nT^n$. Let $M_{\varphi}$ be a
bounded multiplication operator on $L^2(\Sigma)$, then
$q(M_{\varphi})$ is also bounded and
$q(M_{\varphi})=M_{q\circ\varphi}$. By the continuous functional
calculus, for any $f\in C(\sigma(M_{\varphi}))$, we have
$g(M_{\varphi})=M_{g\circ\varphi}$.

\vspace*{0.3cm} {\bf Proposition 2.4.} Let $S=S(E(|u|^2))$ and
$G=S(E(|w|^2))$. If  $f\in C(\sigma(M_{E(|u|^2)}))$ and $g\in
C(\sigma(M_{E(|w|^2)}))$, Then
$$f(T^{\ast}T)=f(0)I+M_{(E(|u|^2))^{-1}.\chi_{S}}\left(M_{f\circ(E(|u|^2)E(|w|2))}-f(0)I\right)M_{\bar{u}}EM_{u}$$
and
$$g(TT^{\ast})=g(0)I+M_{(E(|w|^2))^{-1}.\chi_{G}}\left(M_{g\circ(E(|u|^2)E(|w|2))}-g(0)I\right)M_{w}EM_{\bar{w}}.$$

\vspace*{0.3cm} {\bf Proof.} For all $f\in L^2(\Sigma)$,
$T^{\ast}T(f)=\bar{u}E(|w|^2)E(uf)$ and
$TT^{\ast}(f)=wE(|u|^2)E(\bar{w}f)$. By induction, for each $n\in
\mathbb{N}$,
$$(T^{\ast}T)^n(f)=\bar{u}(E(|w|^2))^n(E(|u|^2))^{n-1}E(uf), \
\ \ (TT^{\ast})^n(f)=w(E(|u|^2))^n(E(|w|^2))^{n-1}E(\bar{w}f).$$
So
$$q(T^{\ast}T)=q(0)I+M_{(E(|u|^2))^{-1}.\chi_{S}}\left(M_{q\circ(E(|u|^2)E(|w|2))}-q(0)I\right)M_{\bar{u}}EM_{u}$$
and
$$q(TT^{\ast})=q(0)I+M_{(E(|w|^2))^{-1}.\chi_{G}}\left(M_{q\circ(E(|u|^2)E(|w|2))}-q(0)I\right)M_{w}EM_{\bar{w}}.$$
By the Weierstrass approximation theorem we conclude that, for
every $f\in C(\sigma(M_{E(|u|^2)}))$ and $g\in
C(\sigma(M_{E(|w|^2)}))$,
$$f(T^{\ast}T)=f(0)I+M_{(E(|u|^2))^{-1}.\chi_{S}}\left(M_{f\circ(E(|u|^2)E(|w|2))}-f(0)I\right)M_{\bar{u}}EM_{u}$$
and
$$g(TT^{\ast})=g(0)I+M_{(E(|w|^2))^{-1}.\chi_{G}}\left(M_{g\circ(E(|u|^2)E(|w|2))}-g(0)I\right)M_{w}EM_{\bar{w}}.$$

\vspace*{0.3cm} {\bf Theorem 2.5.} The unique polar decomposition
of $T=M_wEM_u$ is $U|T|$, where

$$|T|(f)=\left(\frac{E(|w|^{2})}{E(|u|^{2})}\right)^{\frac{1}{2}}\chi_{S}\bar{u}E(uf), \ \ \ \ \ U(f)=\left(\frac{\chi_{S\cap
 G}}{E(|w|^{2})E(|u|^{2})}\right)^{\frac{1}{2}}wE(uf),$$
for all $f\in L^{2}(\Sigma)$.

\vspace*{0.3cm} {\bf Proof.} By Proposition 2.4 we have
$$|T|(f)=(T^{\ast}T)^{\frac{1}{2}}(f)=\left(\frac{E(|w|^{2})}{E(|u|^{2})}\right)^{\frac{1}{2}}\chi_{_S}\bar{u}E(uf).$$
Define  a linear operator $U$ whose action is given by
$$U(f)=\left(\frac{\chi_{S\cap
 G}}{E(|w|^{2})E(|u|^{2})}\right)^{\frac{1}{2}}wE(uf), \ \ \ \ f\in L^2(\Sigma).$$
Then $T=U|T|$ and by Theorem 2.3, $U$ is a partial isometry. Also,
it is easy to see that $\mathcal{N}(T)=\mathcal{N}(U)$. Since for
all $f\in
 L^{2}(\Sigma)$,  $\|Tf\|_2=\|\ |T|f\|_2$, hence $\mathcal{N}(|T|)=\mathcal{N}(U)$ and so
this decomposition is unique.
 $\hfill\Box$

\vspace*{0.3cm} {\bf Theorem 2.6.} The Aluthge transformation of
$T=M_wEM_u$ is
$$\widehat{T}(f)=\frac{\chi_{S}E(uw)}{E(|u|^{2})}\bar{u}E(uf), \ \ \ \ \ \  \ \ \ \  \ \  \ f\in L^{2}(\Sigma).$$

\vspace*{0.3cm} {\bf Proof.} Define operator $V$ on
$L^{2}(\Sigma)$ as
$$Vf=\left(\frac{E(|w|^{2})}{(E(|u|^{2}))^{3}}\right)^{\frac{1}{4}}\chi_{S}\bar{u}E(uf), \ \ \ \ f\in L^{2}(\Sigma).$$
Then we have $V^{2}=|T|$ and so by direct computation we obtain
$$\widehat{T}(f)=|T|^{\frac{1}{2}}U|T|^{\frac{1}{2}}(f)=\frac{\chi_{S}E(uw)}{E(|u|^{2})}\bar{u}E(uf).$$
$\hfill\Box$

\section{ \sc\bf Spectral decomposition }
The normal operators form one of the best understood and most
tractable of classes of operators. The principal reason for this
is the spectral theorem, a powerful structure theorem that answers
many (not all) questions about these operators. In this section we
explore spectral measure and spectral decomposition corresponding
to a normal weighted conditional expectation operator $EM_u$ on
$L^2(\Sigma)$.

\vspace*{0.3cm} {\bf Definition 3.1.} If $X$ is a set, $\Sigma$ a
$\sigma$-algebra of subsets of $X$ and $H$ a Hilbert space, a
spectral measure for $(X, \Sigma, H)$ is a function
$\mathcal{E}:\Sigma\rightarrow B(H)$ having
the following properties.\\
(a) $\mathcal{E}(S)$ is a projection.\\
(b) $\mathcal{E}(\emptyset)=0$ and $\mathcal{E}(X)=I$.\\
(c) If $S_1, S_2\in \Sigma$. $\mathcal{E}(S_1\cap S_2)=\mathcal{E}(S_1)\mathcal{E}(S_2)$.\\
(d) If $\{S_n\}^{\infty}_{n=0}$ is a sequence of pairwise disjoint
sets in $\Sigma$, then
$$\mathcal{E}(\cup^{\infty}_{n=0}S_n)=\Sigma^{\infty}_{n=0}\mathcal{E}(S_n).$$

The spectral theorem says that: For every normal operator $T$ on a
Hilbert space $H$, there is a unique spectral measure
$\mathcal{E}$ relative to $(\sigma(T), H)$ such that
$T=\int_{\sigma(T)}zd\mathcal{E}$, where $z$ is the inclusion map
of $\sigma(T)$ in $\mathbb{C}$.

 J. Herron showed that $\sigma(EM_u)=ess \
range(Eu)\cup\{0\}$, \cite{her}. Also, He has proved that:
 $EM_u$ is normal if and only if $u\in
L^{\infty}(\mathcal{A})$. If $T=EM_u$ is normal, then
$T^n=M_{u^n}E$ and $(T^{\ast})^n=M_{\bar{u}^n}E$. So
$(T^{\ast})^nT^m=M_{(\bar{u})^nu^m}E$ and
$$P(T,T^{\ast})=\sum^{N,M}_{n,m=0}\alpha_{n,m}T^m(T^{\ast})^n=\sum^{N,M}_{n,m=0}\alpha_{n,m}\bar{u}^nu^mE=P(u,\bar{u})E=EP(u,\bar{u}).$$
Where $p(z,t)=\sum^{N,M}_{n,m=0}\alpha_{n,m}z^mt^n$. If
$q(z)=\sum^{N}_{n=0}\alpha_{n}z^n$, then
$q(T)=\sum^{N}_{n=0}\alpha_{n}u^nE$. Hence by the Weierstrass
approximation theorem  we have $f(T)=M_{f(u)}E$, for all $f\in
C(\sigma(EM_u))$. Thus $\phi:C(\sigma(EM_u))\rightarrow
C^*(EM_u,I)$, by $\phi(f)=M_{f(u)}E$, is a unital
$\ast$-homomorphism. Moreover, by Theorem 2.1.13 of \cite{mur},
$\phi$ is also a uniqe $\ast$-isomorphism such that
$\phi(z)=EM_u$, where
$z:\sigma(EM_u)\rightarrow \mathbb{C}$ is the inclusion map.\\

 If $EM_u$ is normal and compact, then
$\sigma(EM_u)=\{0\}\cup\{\lambda_n\}_{n\in \mathbb{N}}$ where
$\lambda_n\neq0$ for all $n\in \mathbb{N}$. So, for each $n\in
\mathbb{N}$
\begin{align*}
E_{n}&=\{0\neq f\in L^2(\Sigma):E(uf)=\lambda_{n}f\}\\
&=\{0\neq f\in L^2(\mathcal{A}):uf=\lambda_{n}f\}\\
 &=L^2(A_n, \mathcal{A}_n, \mu_n),
\end{align*}
 where $A_n=\{x\in X: u(x)=\lambda_n\}$, $E_0=\{f\in
L^2(\Sigma):E(uf)=0\}$, $\mathcal{A}_n=\{A_n\cap B: B\in
\mathcal{A}\}$ and $\mu_n\equiv\mu\mid_{\mathcal{A}_n}$. It is
clear that for all $n,m \in \mathbb{N}\cup\{0\}$, $E_n\cap
E_m=\emptyset$. This implies that the spectral decomposition of
$EM_u$ is as follows:
$$EM_u=\sum^{\infty}_{n=0}\lambda_{n}P_{E_n},$$
where $P_{E_n}$ is the orthogonal projection onto $E_n$.
 Since $EM_u$ is normal, then
$\sigma(EM_u)=ess \ range(u)\cup\{0\}$. So
$\{\lambda_n\}^{\infty}_{n=0}$ is a resolution of the identity on
$X$.\\

Suppose that $W=\{u\in L^0(\Sigma):E(|u|^2)\in
L^{\infty}(\mathcal{A})\}$. If we set
$\|u\|=\|(E(|u|^2))^{\frac{1}{2}}\|_{\infty}$, then $W$ is a
complete ${\ast}$-subalgebra of $L^{\infty}(\Sigma)$.

In the sequel we assume that, $\varphi:X\rightarrow X$ is
nonsingular transformation i.e, the measure $\mu\circ\varphi^{-1}$
is absolutely continuous with respect to the measure $\mu$, and
$\varphi^{-1}(\Sigma)$ is a sub-$\sigma$-finite algebra of
$\Sigma$. Put $h=d\mu\circ\varphi^{-1}/d\mu$ and
$E^{\varphi}=E^{\varphi^{-1}(\Sigma)}$. For $S\in \Sigma$, let
$\mathcal{E}(S):L^2(\Sigma)\rightarrow L^2(\Sigma)$ be defined by
$$\mathcal{E}(S)(f)=E^{\varphi}M_{\chi_{\varphi^{-1}(S)}}(f),$$ i.e,
$$\mathcal{E}(S)=E^{\varphi}M_{\chi_{\varphi^{-1}(S)}}.$$ $\mathcal{E}$ defines a
spectral measure for $(X, \Sigma, L^2(\mu))$. If $E^{\varphi}M_u$
is normal on $L^2(\Sigma)$, then by Theorem 2.5.5 of \cite{mur},
$\mathcal{E}$ is the unique spectral measure corresponding to
$\ast$-homomorphism $\phi$ that is defined as follows:

$$\phi:C(\sigma(E^{\varphi}M_u))\rightarrow C^*(EM_u,I), \ \ \ \
\phi(f)=E^{\varphi}M_{f(u)}.$$ So, for all $f\in
C(\sigma(E^{\varphi}M_u))$ we have
$$\phi(f)=\int_{X}fd\mathcal{E}.$$
In \cite{cl} it is explored that which sub-$\sigma$-algebras of
$\Sigma$ are of the form $\varphi^{-1}(\Sigma)$ for some
nonsingular transformation $\varphi:X\rightarrow X$.
 These observations establish the following theorem.

\vspace*{0.3cm} {\bf Theorem 3.2.}
 Let $(X, \Sigma, \mu)$ be a $\sigma$-finite measure
space,  $\varphi:X\rightarrow X$ be a nonsingular transformation
and let $u$ be in $L^{\infty}(\varphi^{-1}(\Sigma))$. Consider the
operator $E^{\varphi}M_u$ on $L^2(\Sigma)$. Then the set function
$\mathcal{E}$ that is defined as:
$\mathcal{E}(S)=E^{\varphi}M_{\chi_{\varphi^{-1}(S)}}$ for $S\in
\Sigma$, is a spectral measure. Also, $\mathcal{E}$ has compact
support and
$$E^{\varphi}M_u=\int zd\mathcal{E}.$$

\end{document}